\documentclass[11pt]{article}
\usepackage{graphicx} 
\usepackage[a4paper,
            bindingoffset=0.2in,
            left=3cm,
            right=3cm,
            bottom=3cm,
            top = 3.3cm]{geometry}
\usepackage[toc]{appendix}
\usepackage{amsthm}
\usepackage{amsmath}
\usepackage{amsfonts}
\usepackage{array}
\usepackage{listings}
\usepackage[maxbibnames=99]{biblatex} 
\usepackage{hyperref}
\hypersetup{
    pdftitle={Implementing Hadamard Matrices in SageMath},
    colorlinks=true,
    linkcolor=black,
    urlcolor=blue,
    citecolor=black,
    }

\newtheorem*{definition}{Definition}
\newtheorem{theorem}{Theorem}

\title{Implementing Hadamard Matrices in SageMath}
\date{}

\begin{document}

\begin{titlepage}
    \centering
    \begin{figure}
        \centering
        \includegraphics[width=200px]{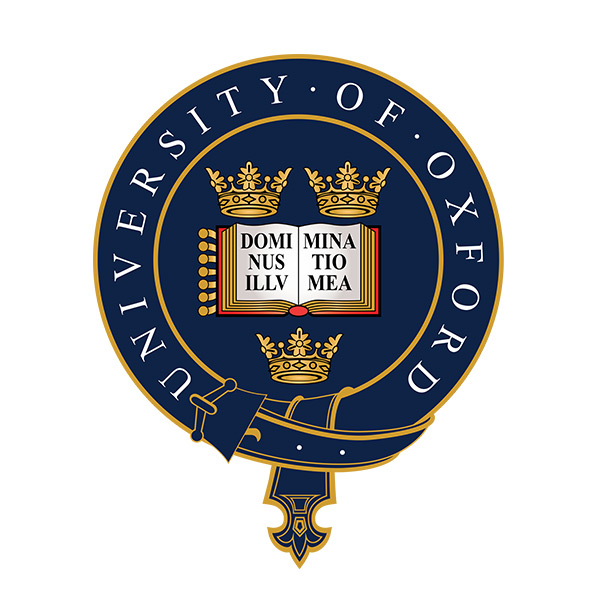}
    \end{figure}
    
    \Huge
    \hrule
    \vspace*{15px}
    \textbf{Implementing Hadamard Matrices in SageMath}
    \vspace*{10px}
    \hrule
    \vspace*{1cm}

    \Large
    Matteo Cati \\
    \vspace*{10px}
    Keble College \\
    University of Oxford 
    
    \vspace*{15px}
    \textit{Supervisor:} Dr Dmitrii Pasechnik 
    
    \Large
    \vspace*{100px}
    \textit{Honour School of Computer Science, Part B}\\
    \vspace*{15px}
    Trinity 2023

    \vspace*{\fill}
\end{titlepage}
\clearpage

\pagenumbering{arabic}
\section*{Abstract}

Hadamard matrices are $(-1, +1)$  square matrices with mutually orthogonal
rows. The Hadamard conjecture states that Hadamard matrices of order $n$ exist
whenever $n$ is $1$, $2$, or a multiple of $4$. However, no construction is known
that works for all values of $n$, and for some orders no Hadamard matrix has
yet been found. Given the many practical applications of these matrices, it
would be useful to have a way to easily check if a construction for a Hadamard
matrix of order $n$ exists, and in case to create it. This project aimed to
address this, by implementing constructions of Hadamard and skew Hadamard
matrices to cover all known orders less than or equal to $1000$ in SageMath, an
open-source mathematical software. Furthermore, we implemented some additional
mathematical objects, such as complementary difference sets and T-sequences,
which were not present in SageMath but are needed to construct Hadamard
matrices. 

This also allows to verify the correctness of the results given in the
literature; within the $n\leq 1000$ range, just one order, $292$, of a skew Hadamard
matrix claimed to have a known construction, required a fix.

\clearpage
\tableofcontents
\clearpage
\section{Introduction}
\label{intro-chapter}
Hadamard matrices were defined for the first time by Sylvester in 1867. Since then, they have been used in many practical applications, such as quantum computing, image analysis, and signal processing.

However, no construction is currently known that works for matrices of every size. Therefore, whenever a Hadamard matrix of a specific order is needed, one must try multiple constructions until the correct one is  found (if it is known).

This project aims to alleviate this problem, by providing all the necessary constructions in SageMath, an open-source mathematical software written in Python. 
Although a few constructions were already present in SageMath (see chapter \ref{background-chapter}), they provided matrices for infinite series which only covered a few orders. We added many new constructions so that Hadamard and skew Hadamard matrices of all known orders less than or equal to 1000 are now available. 

The code can be found on the SageMath GitHub repository\footnote{\url{https://github.com/sagemath/sage}}. In particular, the functions that I have written have been added to the repository with multiple separate Pull Requests\textsuperscript{2-7}. This was done to follow the development guidelines of SageMath.
\footnotetext[2]{\url{https://github.com/sagemath/sage/issues/32267}}
\footnotetext[3]{\url{https://github.com/sagemath/sage/issues/34690}}
\footnotetext[4]{\url{https://github.com/sagemath/sage/issues/34807}}
\footnotetext[5]{\url{https://github.com/sagemath/sage/pull/34985}}
\footnotetext[6]{\url{https://github.com/sagemath/sage/pull/35059}}
\footnotetext[7]{\url{https://github.com/sagemath/sage/pull/35211}}

There were a few  particularly important aspects to take into consideration. First of all, given that new constructions are found fairly often, the functions provided should be easily extendable. Furthermore, we must be able to guarantee the correctness of the constructions implemented. Lastly, given that the code is part of open-source software, it should contain clear documentation and must adhere to the coding style adopted by the SageMath community.

\subsection{Structure of the Report}
This report is divided into the following chapters:
\begin{itemize}
    \item Chapter \ref{background-chapter} contains a definition of Hadamard matrices and a brief description of some methods that were already implemented in SageMath;
    \item Chapter \ref{design-chapter} gives a high-level description of the functions implemented;
    \item Chapter \ref{NA-chapter} describes the construction of some sequences with zero nonperiodic autocorrelation;
    \item Chapter \ref{diff-sets-chapter} gives some constructions of difference sets;
    \item Chapter \ref{hadamard-chapter} contains a description of the constructions for Hadamard (and skew Hadamard) matrices that have been implemented;
    \item Chapter \ref{conclusion-chapter} contains some final considerations.
\end{itemize}

\clearpage
\section{Background}
\label{background-chapter}
A $n \times n$ matrix is called a \textit{Hadamard matrix} if its entries are all $-1$, $+1$ and the rows are mutually orthogonal. Equivalently, they are $(-1, +1)$ matrices which satisfy the equation:
\begin{displaymath}
    HH^\top = H^\top H = nI
\end{displaymath}
In particular, the latter equality implies that $|det(H)| = n^{n/2}$, i.e. the determinant is maximal.
Matrices of this type were described for the first time by Sylvester \cite{Sylvester1867}, and were studied further by Jacques Hadamard \cite{Hadamard1893}.

Additionally, a Hadamard matrix is  \textit{skew} if it has the form $H=S+I$, where $I$ is the identity and $S$ is skew-symmetric: $S^\top = -S$.

The Hadamard conjecture states that a Hadamard matrix of order $n$ exists if and only if $n=1, 2$ or a multiple of four. It is easy to see that for a Hadamard matrix of order $n > 2$ to exist, $n$ must be a multiple of four. However, the converse has not yet been proven. Currently, Hadamard matrices are known for most orders less than 1000, with the only exceptions being 668, 716, and 892.

A similar conjecture has been proposed for skew Hadamard matrices. Again, skew Hadamard matrices are not known for all orders, and in particular no construction is known for the following values of $n < 1000$ (see \cite{Djokovic2016}):
\begin{align*}
   \centerline{356, 404, 428, 476, 596, 612, 668, 708, 712, 716, 764, 772, 804,} \\
   \centerline{808, 820, 836, 856, 892, 900, 916, 932, 940, 952, 980, 996}
\end{align*}
 
In the following sections, we will see some constructions of  Hadamard matrices that had already been implemented in SageMath.

\subsection{Paley's Constructions}
In 1933, Paley \cite{Paley1933} discovered two constructions of Hadamard matrices. They are as follows:

\begin{theorem}
\label{paleyI-theo}
Let $q$ be a prime power, with $q \equiv 3 \mod 4$. Then there is a skew Hadamard matrix of order $q+1$. \qed
\end{theorem}

\begin{theorem}
Let $m = 2(q+1)$, where $q \equiv 1 \mod 4$ is a prime power. Then there is a Hadamard matrix of order $m$. \qed
\end{theorem}

\subsection{Doubling Construction}
Sylvester \cite{Sylvester1867} proved that if $n$ is the order of a Hadamard matrix, then there exists a Hadamard matrix of order $2^tn$ for all values of $t\ge0$.

\begin{theorem}
Let $H$ be a Hadamard matrix of order $n$. Then the matrix $H' = \begin{pmatrix}
H & H\\
H & -H
\end{pmatrix}$ is an Hadamard matrix of order $2n$. \qed
\end{theorem}

Furthermore, Seberry \cite{Wallis1971} described a similar construction for skew Hadamard matrices:
\begin{theorem}
    Suppose $H_n = S + I_n$ is a skew Hadamard matrix of order $n$. Then

    \[
    H_{2n} = 
    \begin{pmatrix}
    S+I_n & S+I_n\\
    S-I_n & -S+I_n
    \end{pmatrix}
    \]

    is a skew Hadamard matrix of order $2n$. \qed
\end{theorem}

These theorems are particularly important, because they imply that when creating new Hadamard matrices of order $4n$ we only need to worry about the cases when $n$ is odd.

\subsection{Regular Symmetric Hadamard Matrices with Constant Diagonal}
Regular symmetric Hadamard matrices with constant diagonal (or RSHCD) are symmetric $n \times n$ Hadamard matrices such that:
\begin{itemize}
    \item All values on the main diagonal are equal to a value  $\delta \in \{-1, +1\}$;
    \item All rows sums are equal to $\delta \epsilon \sqrt{n}$, with $\epsilon \in \{-1, +1\}$.
\end{itemize}

All known RSHCD of order $n\le 1000$ had already been implemented in SageMath \cite{Cohen2017}. However, they  only cover the orders:
\begin{align*}
    4, 16, 36, 64, 100, 144, 196, 256, 324, 400, 576, 676, 784, 900
\end{align*}

\clearpage
\section{Program Design}
\label{design-chapter}
All the functions implemented have a common structure. First of all, given that the code is part of an open-source software, it was of the utmost importance that good documentation was provided. In SageMath, documentation is written using the functions' docstrings. These contain a general description of the implementation, with some references to the relevant papers. Subsequently, one can find a list of the input parameters, together with a short explanation of each one, and a description of the resulting output. Lastly, the documentation provides a few examples of how to use the function.

\subsection{Testing}
In order to look for possible bugs in the code, the docstrings of every function contain some unit tests, which are run by GitHub Continuous Integration tools whenever new code is added to the repository. On average, every function I wrote contains between five and ten tests, depending on its complexity. 

Writing these tests is particularly important because, given the large size of the SageMath repository, some new pieces of code may introduce bugs in seemingly unrelated parts of the software.

\subsection{Common Function Parameters}
Since SageMath is a mathematical software, the correctness of the implementation must be guaranteed. Although, as seen in the previous section, many tests have been written for each function, an additional precaution was taken. 

In particular, every function that creates a mathematical object runs a check on the object created, to confirm that the result is correct.
However, such checks often have quadratic (or worse) complexity, which may add considerable overhead when the input is large. Therefore, the function parameters contain a boolean flag, and the check will be skipped if the parameter is set to false by the user.

Furthermore, it is often not easy to see if a construction can be applied to some input or not (for example, it may depend on whether data for such construction is contained in SageMath). To make the code more user-friendly, in all such cases a boolean parameter called {\tt existence} is present. When this parameter is set to true, the function will not compute the requested object: instead, it returns true if the function can be applied to the given input, and false otherwise. 

As an example, Fig. \ref{ex-parameters-code} shows the pseudocode of a generic function computing Hadamard matrices.

\begin{lstlisting}[float, label=ex-parameters-code, caption={Functions structure}]
def hadamard_matrix(n, existence=False, check=True):
    if existence:
        if Hadamard matrix of order n can be created:
            return True
        return False

    if Hadamard matrix of order n cannot be created:
        raise ValueError()

    H := create Hadamard matrix of order n
    if check:
        raise exception if H is not a Hadamard matrix
    return H
    
\end{lstlisting}

\subsection{Utility Functions}
\label{utility-funcs-section}
In the previous section, we have seen that most functions need a way to check if the result is the correct mathematical object. To avoid code duplication, I have implemented some utility functions that execute this check. For example, the function {\tt are\_complementary\_difference\_sets} is used to check if the given sets are complementary difference sets.

Furthermore, these functions contain a boolean parameter { \tt verbose}. When this is true, the function will print to the console a string explaining why the result is ``false". This is particularly useful when we need to check that an object satisfies multiple properties: for example, if we want to check if a matrix is a skew Hadamard matrix and the function {\tt is\_skew\_hadamard\_matrix} returns false, setting {\tt verbose} to true will tell us if the matrix is not Hadamard, or if it is Hadamard but not skew.

Lastly, many of the mathematical objects detailed in this project can be obtained from multiple constructions. In these cases, the code contains a utility function (e.g. {\tt hadamard\_matrix}) that tries to find a result by checking each of the given constructions in order. 

\clearpage
\section{Sequences with Zero Nonperiodic Autocorrelation}
\label{NA-chapter}
Of particular interest for the construction of Hadamard matrices are complementary sequences. 

In general, given a sequence $A=(a_1, a_2, ..., a_n)$ the nonperiodic autocorrelation function $N_A(j)$ is defined by \cite{Kharaghani2005}:
\begin{align*}
    N_A(j) &= \sum_{i=1}^{n-j}a_ia_{i+j} &\text{ for } 0 \le j \le n \\
    N_A(j) &= 0 &\text{ for } j > n 
\end{align*}

According to \cite{Seberry2017}, a family  $X=\{A_1, A_2, ..., A_k\}$ of integer sequences of length $n$ is complementary if it satisfies:
\begin{align*}
    \sum_{i=1}^{k}N_{A_i}(j) = 0 \quad j\in \{1, 2, ..., n-1\}
\end{align*}

\subsection{Turyn Sequences}
\label{TurynSeqs}
The first sequences that I implemented are Turyn sequences. These are families composed by four $(-1, +1)$ complementary sequences of length $l, l, l-1, l-1$, with the form (see Definition 7.4 of \cite{Seberry2017}):
\begin{align*}
    X &= \{x_1 = 1, x_2, x_3, ..., x_{l-1}, x_l = -1\}\\
    U &= \{u_1, u_2, ..., u_{l-1}, u_l = 1\}\\
    Y &= \{y_1, y_2, .., y_{l-1}\}\\
    V &= \{v_1,v_2, ..., v_{l-1}\}
\end{align*}

These sequences, which were constructed for the first time by Turyn in \cite{Turyn1974},  are known for $l = 2, 3, 4, 5, 6, 7, 8, 13, 15$.

\subsection{Turyn Type Sequences}
Some other sequences used for constructing Hadamard matrices are Turyn type sequences. Formally \cite{Kharaghani2005},
\begin{definition}
Four $(-1, +1)$ sequences $X, Y, Z, W$ of lengths $n, n, n, n-1$ are said to be of
Turyn type if
\begin{displaymath}
    N_X(j) + N_Y(j) + 2N_Z(j) + 2N_W(j) = 0 \text{ for }  j\ge 1
\end{displaymath}
\end{definition}

The sequences of this type that I implemented in SageMath are listed in  \cite{TurynTypeSeqs, Kharaghani2005}.

\subsection{Base Sequences}
The Turyn sequences described in section \ref{TurynSeqs} can be generalised into the concept of base sequences \cite{Kharaghani2005}:
\begin{definition}
    Four $(-1, 1)$ sequences $A, B, C, D$ of lengths $n + p, n + p, n, n$ are called base sequences if
    \begin{displaymath}
        N_A(j) + N_B(j) + N_C(j) + N_D(j) = 0 \text{ for }  j\ge 1
    \end{displaymath}
\end{definition}

Clearly, Turyn sequences can be seen as base sequences with $p=1$.
Additionally, base sequences of lengths $2n-1, 2n-1, n, n$ can be constructed from Turyn type sequences of lengths $n,n,n,n-1$ \cite{Kharaghani2005}:
\begin{theorem}
    If $X, Y, Z, W$ are Turyn type sequences of lengths $n, n, n, n - 1$, then the sequences $A = Z;W$, $B = Z;-W$, $C = X$, $D = Y$ are base sequences of lengths $2n - 1, 2n - 1, n, n$. \qed
\end{theorem}

Here the notation $A; B$ is used to mean sequence $A$ followed by sequence $B$.

\subsection{T-sequences}
Lastly, I used the sequences defined so far to implement some constructions of T-sequences. They are defined by \cite{Kharaghani2005} as:
\begin{definition}
    Four $(-1, 0, 1)$ sequences $A, B, C, D$ of length $n$ are called T-sequences if

\begin{displaymath}
    N_A(j) + N_B(j) + N_C(j) + N_D(j) = 0  \text{\: for \:}  j\ge 1
\end{displaymath}
    and in each position exactly one of the entries of $A, B, C, D$ is nonzero.
\end{definition}

The first way to generate T-sequences is by using base sequences. From \cite{Kharaghani2005}:

\begin{theorem}
Suppose $A, B, C, D$ are base sequences of lengths $n+p, n+p, n, n$. Then, the following are T-sequences of length $2n+p$:
\begin{align*}
    T_1 &= \frac{1}{2}(A+B); 0_n \\
    T_2 &= \frac{1}{2}(A-B); 0_n \\
    T_3 &= 0_{n+p};\frac{1}{2}(C+D)\\
    T_4 &= 0_{n+p}; \frac{1}{2}(C-D)
\end{align*}\qed
\end{theorem}
The notation used is as follows. Suppose $A = (a_1, a_2, .., a_k)$, $B = (b_1, b_2, ..., b_k)$ are two sequences, and $c$ is a scalar. Then:

\begin{itemize}
    \item $c_n$ is a sequence of length $n$ whose elements are all $c$;
    \item $A;B = (a_1, ..., a_k, b_1, ..., b_k)$;
    \item $A+B = (a_1+b_1, a_2+b_2, ..., a_k+b_k)$, and similarly $A-B = (a_1-b_1, a_2-b_2, ..., a_k-b_k)$;
    \item $cA = (ca_1, ca_2, ..., ca_k)$.
\end{itemize}

In addition, T-sequences of length $4n-1$ can be constructed from Turyn sequences of lengths $n, n, n-1, n-1$. From Theorem 7.7 of \cite{Seberry2017}:

\begin{theorem}
\label{theo-constrTsequences-4n-1}
    Suppose $A, B, C, D$ are four Turyn sequences of lengths $n, n, n-1, n-1$. Then, the following are T-sequences of length $4n-1$:
    \begin{align*}
        T_1 &= 1; 0_{4n-2} \\
        T_2 &= 0; A/C; 0_{2n-1} \\
        T_3 &= 0_{2n}; B/0_{n-1} \\
        T_4 &= 0_{2n}; 0_n/D
    \end{align*}\qed
\end{theorem}

The notation used in Theorem \ref{theo-constrTsequences-4n-1} is the same as for the previous construction, with the addition that if we have two sequences $A = (a_1, a_2, ..., a_k)$ and $B = (b_1, b_2, ..., b_{k-1})$ then $A/B = (a_1, b_1, a_2, b_2, ..., a_{k-1}, b_{k-1}, a_k)$.

\clearpage
\section{Difference Sets}
\label{diff-sets-chapter}
Another concept often used in the construction of Hadamard matrices is that of difference sets.

\begin{definition}
    A $(v, k, \lambda)$ difference set is a subset $D$ of size $k$ of an additive group $G$ of order $v$ such that every nonidentity element of $G$ can be expressed as a difference $d_1 - d_2$ of elements of $D$ in exactly $\lambda$ ways.
\end{definition}

We will focus on the construction of some variations of difference sets: relative difference sets, supplementary difference sets, and complementary difference sets.

\subsection{\textit{m}-sequences}
\label{m-sequences-sect}
Given a Galois field $G$ of order $q = p^m$, an \textit{m}-sequence is a periodic (infinite) sequence with period $q^n-1$, with elements from $G$.

These sequences, which will be used in the construction of relative difference sets (section \ref{rel-diff-set}), can be created as described by Zierler in \cite{Zierler1959}. In particular, given $n+1$ values $(c_0, c_1, ..., c_n)$, and the first $n$ values $a_1, a_2, ..., a_n$ we define the entire sequence as follows:

\begin{displaymath}
    a_k = -c_0^{-1}\sum_{i = 1}^{n}c_ia_{k-i} \text{ for } k > n
\end{displaymath}

Furthermore, when $f(x) = c_0 + c_1x + ... + c_n x^n$ is a primitive polynomial, Zierler proved that there are values $a_1, ..., a_n$ which construct an \textit{m}-sequence.

In fact, \cite{Mitra2008} showed that when $q=2$, we can use as initial sequence $(1, 0, 0, ..., 0)$. This can be generalised to any value of $q$. We have therefore a complete algorithm for finding an \textit{m}-sequence (Fig. \ref{mSequences-code}). 

This algorithm assumes that there exists a function {\tt find\_primitive\_poly} for finding a primitive polynomial over a given ring. This was not available in SageMath, so I implemented it by repeatedly creating a random irreducible polynomial until one that is also primitive is found (see Fig. \ref{primitivePoly-code}). Although no guarantee of the efficiency of this function can be given, it has proven to be fast in practice.
 
\begin{lstlisting}[caption={Function constructing \textit{m}-sequences}, label=mSequences-code, float]
def create_m_sequence(q, n):
    K = GF(q)
    T = PolynomialRing(K, 'x')
    primitive = find_primitive_poly(T, n)
    
    coeffs = primitive.coefficients()
    exps = primitive.exponents()

    seq = [1] + [0]*(n-1)
    while len(seq) < q**n - 1:
        nxt = 0
        for i, coeff in zip(exps[1:], coeffs[1:]):
            nxt += coeff * seq[-i]
        seq.append(-coeffs[0].inverse() * nxt)
    return seq
\end{lstlisting}

\begin{lstlisting}[caption={Function creating primitive polynomials}, label=primitivePoly-code, float]
def find_primitive_poly(T, n):
    primitive = T.irreducible_element(n, algorithm='random')
    while not primitive.is_primitive():
        primitive = T.irreducible_element(n, algorithm='random')
    return primitive
\end{lstlisting}

\subsection{Relative Difference Sets}
\label{rel-diff-set}
Elliot and Butson \cite{Elliott1966} gave the following definition of relative difference sets:
\begin{definition}
    A set $R$ of elements in a group $G$ of order $mn$ is a difference set of $G$ relative to a normal subgroup $H$ of order $n \not = mn$ if the collection of differences $r - s$, $r, s \in R$, $r \not = s$ contains only elements of $G$ which are not in $H$, and contains every such element exactly $d$ times.
\end{definition}
Such sets are denoted as $R(m, n, k, d)$. Furthermore, we call the set cyclic or abelian if  the group $G$ has the corresponding property. 

A first construction for relative difference sets which uses \textit{m}-sequences is described in \cite{Elliott1966}:
\begin{theorem}
\label{relDiffSet-m-seq-theo}
    For each \textit{m}-sequence over a field of $q=p^m$ elements, there exists a cyclic $R((q^N-1)/(q-1), q-1, q^{N-1}, q^{N-2})$ where $q^N-1$ is the period of the \textit{m}-sequence. \qed
\end{theorem}

Given a \textit{m}-sequence $a = \{a_i\}_{i=0}^\infty$ of period $q^N-1$, the relative difference set over the group of integers modulo $q^N-1$ is:
\begin{displaymath}
    R = \{i; 0\le i < q^N-1 | a_i = 1\}
\end{displaymath}

Using a second construction, we can create relative differences sets with parameters $R((q^N-1)/(q-1), n, q^{N-1}, q^{N-2}d)$, where $q$ is a prime power and $nd = q-1$.

This is possible because of a theorem from \cite{Elliott1966}: 
\begin{theorem}
    If R is an $R(m, n, k, d)$ and if $\sigma$ is a homomorphism of $G$ onto $\sigma(G)$ with kernel $K \subset H$, then $\sigma(R)$ is an $R(m, s, k, td)$ of $\sigma(G)$ relative to $\sigma(H)$, where $n=ts$, and $t$ is the order of $K$. \qed
\end{theorem}

First, we create a $R((q^N-1)/(q-1), q-1, q^{N-1}, q^{N-2})$ relative difference set $S_1$ using the previous construction. Then, we create a homomorphism $\sigma$ whose kernel has order $d$. Finally, by applying this homomorphism to $S_1$ we get the desired relative difference set.

\begin{lstlisting}[caption={Function computing $R((q^N-1)/(q-1), q-1, q^{N-1}, q^N-2)$}, label=relDiffSet-homo-code, float]
def relative_difference_set_from_homomorphism(q, N, d):
    G := AdditiveAbelianGroup of order q^N - 1
    K := subgroup of G of order d

    sigma = homomorphism with kernel K
    
    diff_set := relative_difference_set_from_m_sequence(q, N)
    second_diff_set := [sigma(x) for x in diff_set]
    return second_diff_set
\end{lstlisting}

Note that if $D$ is a relative difference set over a group $G$, then for any value $t \in G$ the set $\{t+d\text{ } | \text{ }d \in D\}$ is also a relative difference set \cite{Spence1975} and is called a translate of $D$.

Lastly, we define an additional property of relative difference sets \cite{Spence1975}:
\begin{definition}
    We say that a relative difference set $D$ is fixed by $t\in G$ if 
    \begin{displaymath}
        \{td \text{ } | \text{ } d\in D\} = D
    \end{displaymath}
\end{definition}
\subsection{Supplementary Difference Sets}
Another variation of difference sets are supplementary difference sets. From Definition 4.3 of \cite{Seberry2017}:

\begin{definition}
    Let $S_1, S_2, ..., S_n$ be subsets of $G$, an additive abelian group
of order $v$. Let $|S_i| = k_i$. If the equation $g = r - s$, $r, s \in S_i$ has exactly $\lambda$ solutions for each non-zero element $g$ of $g$, then $S_1, S_2, ..., S_n$ will be called $n - \{v;k_1, k_2, ...,k_n;\lambda\}$ supplementary difference sets (or SDS). If $k_1 = k_2 = ... = k_n = k$, we call it a $n - \{v;k;\lambda\}$ SDS.
\end{definition}

Furthermore, if $S_1 \cap -S_1 = \emptyset$ and $S_1 \cup -S_1 = G\setminus \{0\}$ we say that $S_1$ is skew, and the sets are called skew SDS.

\subsubsection{SDS from Relative Difference Sets}
\label{SDS-from-rel-diff-sets-section}
A construction for an infinite family of $4-\{2v; v, v+1, v, v; 2v\}$ SDS is given by Spence in \cite{Spence1975}. 
\begin{theorem}
    If $q$ is an odd prime power for which there exists an integer $s > 0$ such that $(q - (2^{s+1} + 1))/2^{s+1}$ is an odd prime power, then there exists a $4-\{2v; v, v+1, v, v; 2v\}$ SDS, where $v = (q-1)/2$. \qed
\end{theorem}

The first step of this construction is to use Theorem \ref{relDiffSet-m-seq-theo} (with $N=2$) to obtain a relative difference set $D$ with parameters $R(q+1, q-1, q, 1)$. Then, we compute a translate of $D$ which is fixed by $q$. In fact, Spence showed that such a set always exists. 

Furthermore, he also showed that this set must contain an element congruent to $0 \mod q+1$, and therefore we can get a new relative difference set $D_1$ (fixed by $q$) where such element is $0$ by translating 
the set by a suitable multiple of $q+1$. 

Finally, Spence showed how we can use $D_1$ to construct four polynomials $\psi_1(x)$, $\psi_2(x)$, $\psi_3(x)$, $\psi_4(x)$ with the form:

\begin{align*}
    \psi_k(x) &= \sum_{s \in S_k} x^s \quad 1 \le k \le 4
\end{align*}

The four sets $S_1, S_2, S_3, S_4$ are SDS with parameters $4-\{2v; v, v+1, v, v; 2v\}$.

\subsubsection{Computation of SDS}
\label{sds-computation-section}
Particularly important are also the SDS with parameters $4-\{n; k_1, k_2, k_3, k_4; n-\sum_{i=1}^4k_i\}$.

Although no general construction is known for these sets, they have been computed for many values of $n$  \cite{Djokovic1992a, Djokovic1992b, Djokovic1992c}. In most cases, they are defined as subsets of the group $G$ of residues modulo $n$. In particular, papers that use this notation usually provide a subset $H$ of $G$ with $|H| =k$, a set of values $C = \{c_i \text{ } | \text{ }  0 \le i \le (n-1)/(2k)\}$ and four indices sets $J_1, J_2, J_3, J_4$.

Then, from the set $H$ we construct $(n-1) / k$ subsets of $G$:
\begin{align*}
    \alpha_{2i} &= \{c_i h \text{ }|\text{ } h \in H \} \\
    \alpha_{2i+1} &= \{-x \text{ }|\text{ } x \in \alpha_{2i}\}
\end{align*}
and the SDS $S_1, S_2, S_3, S_4$ are:
\begin{displaymath}
    S_i = \bigcup_{j \in J_i} \alpha_j \quad 1\le i \le 4
\end{displaymath}

I have generalised this algorithm (summarised in Fig. \ref{sds-from-cosets-code}), so that it is possible to add $0$ to some subsets, and it allows to specify some $\alpha_i$ directly as a set.

\begin{lstlisting}[caption={Basic algorithm for computing SDS}, label=sds-from-cosets-code, float]
def construction_sds(n, H, indices, subsets_gen):
    Z = Zmod(n)

    subsets = []
    for el in subsets_gen:
        even_sub = {x*el for x in H}
        odd_sub = {-x for x in even_sub}
        subsets.append(even_sub)
        subsets.append(odd_sub)

    def generate_set(index_set, subs):
        S = set()
        for idx in index_set:
            S = S.union(subs[idx])
        return S

    S1 = generate_set(indices[0], subsets)
    S2 = generate_set(indices[1], subsets)
    S3 = generate_set(indices[2], subsets)
    S4 = generate_set(indices[3], subsets)
    
    return [S1, S2, S3, S4]
\end{lstlisting}

Lastly, I have implemented some additional constructions. In particular, \cite{Djokovic1994a} defines some skew SDS over a group of polynomials, and \cite{Djokovic2016} gives some skew SDS where $S_1$ is the Paley-Todd difference set.

Overall, skew SDS of this type are now available in SageMath for $n$ equal to:
\begin{align*}
    \centerline{37, 39, 43, 49, 65, 67, 73, 81, 93, 97, 103, 109, 113, 121, 127, 129, 133,} \\
    \centerline{145, 151, 157, 163, 169, 181, 213, 217, 219, 239, 241, 247, 267, 331, 631}
\end{align*}

Additional non-skew SDS are also available for $n=191, 251$.

\subsection{Complementary Difference Sets}
\label{comp-diff-set-section}
Finally, we define complementary difference sets \cite{Szekeres1969}:
\begin{definition}
    Two subsets $A$ and $B$ of an additive abelian group $G$ of order $2m+1$ will be called complementary difference sets in $G$ if
\begin{itemize}
    \item $A$ contains $m$ elements;
    \item $\alpha \in A$ implies $-\alpha \not \in A$;
    \item for each $\delta \in G$, $\delta \not = 0$ the equations
    \begin{align*}
        \delta = \alpha_1 - \alpha_2 \\
        \delta = \beta_1 - \beta_2 
    \end{align*}
    have altogether $m-1$ distinct solution vectors.
\end{itemize}
\end{definition}

Using the notation from the previous section, $A, B$ are complementary difference sets if they are skew SDS with parameters $2-\{2m+1; m, m; m-1\}$. I have implemented three different constructions.

The first construction creates complementary difference sets over a group of order $q$, when $q$ is a prime power congruent to $3 \mod 4$. Given a Galois field $G$ of order $q$, Szekeres \cite{Szekeres1971} showed that the two sets $A=B$ which contain the nonzero squares in $G$ are complementary difference sets over the additive group of $G$.

A second function allows to create complementary difference sets over a group of order $q = p^t$, where $p$ is a prime congruent to $5 \mod 8$, and $t \equiv 1, 2, 3 \mod 4$. In \cite{Szekeres1969}, Szekeres gave a construction for the case $t \equiv 1 \mod 2$, which uses a Galois field $G$ of order $q$. 

Let $\rho$ be a generator of the multiplicative group of $G$, and $C_0$ be the set of non-zero fourth powers of $G$. Then, define the sets $C_1, C_2, C_3$ as $C_i = \rho^i C_0$.
The complementary difference sets over the additive group of $G$ are:
\begin{align*}
    A &= C_0 \cup C_1\\
    B &= C_0 \cup C_3
\end{align*}

In 1971, Szekeres published a new paper \cite{Szekeres1971}, which covers the case $t \equiv 2 \mod 4$. Let $G$ be the Galois field of order $q$. Then, if $\rho$ is a multiplicative generator of $G$, let $C_0$ be the set of nonzero eighth powers of $G$. Define $C_i = \rho^i C_0$ for $1 \le i \le 7$; the complementary difference sets are:
\begin{align*}
    A &= C_0 \cup C_1 \cup C_2 \cup C_3 \\
    B &= C_0 \cup C_1 \cup C_6 \cup C_7
\end{align*}

The last function creates complementary difference sets over a group of order $n=2m+1$, when $4m+3$ is a prime power \cite{Szekeres1969}. Let $\rho$ be a primitive element of the Galois field of order $4m+3$. Then, define $Q$ to be the set of squares in the Galois field. The following sets are complementary difference sets in the group of integers modulo $n$:
\begin{align*}
    A &= \{ a \text{ } | \text{ } 0 \le a \le n \land \rho^{2a}-1 \in Q \}\\
    B &= \{ b \text{ } | \text{ } 0 \le b \le n \land -\rho^{2b}-1 \in Q \}\\
\end{align*}

\clearpage
\section{Hadamard Matrices}
\label{hadamard-chapter}
\subsection{Williamson Construction}
At the start of my project, the first order for which no Hadamard matrix was present in SageMath was $116$. The construction for this order is due to Williamson \cite{Hall1988}: he proved that the matrix $H$ given by 

\begin{displaymath}
H = \begin{pmatrix}
A & B & C & D \\
B & A & D & -C \\
C & -D & A & B \\
D & C & -B & A
\end{pmatrix}
\end{displaymath}

is a Hadamard matrix of order $4n$, if $A, B, C, D$ are $n\times n$ symmetric circulant matrices (i.e. matrices where every row contains the same entries as the previous row, but rotated to the right by one element) that commute with each other. Furthermore, they must satisfy the condition:

\begin{displaymath}
    A^2+B^2+C^2+D^2 = 4nI
\end{displaymath}

Such matrices are called Williamson matrices.

\subsection{Goethals-Seidel Array}
Goethals and Seidel \cite{Goethals1970} discovered that given four $(-1, +1)$ matrices $A, B, C, D$ of order $n$ such that
\begin{displaymath}
    AA^\top+BB^\top+CC^\top+DD^\top = 4nI
\end{displaymath}
then a Hadamard matrix of order $4n$ can be constructed by plugging them into the Goethals-Seidel array:
\begin{displaymath}
    GS(A, B, C, D) = \begin{pmatrix}
        A & BR & CR & DR \\
        -BR & A & D^\top R & -C^\top R \\
        -CR & -D^\top R & A & B^\top R \\
        -DR & C^\top R & -B^\top R &  A
    \end{pmatrix}
\end{displaymath}
Here, $R$ is the $n\times n$ permutation matrix with all ones on the anti-diagonal.

Furthermore, if $A = S+I$ with $S$ a skew-symmetric matrix, the Goethals-Seidel array will give a skew Hadamard matrix.

\subsubsection{Hadamard Matrices from Supplementary Difference Sets}
The Goethals-Seidel array has been used extensively to construct Hadamard matrices (both skew and non-skew) of order $4v$ from SDS over groups of size $v$. As explained in \cite{Djokovic2008}, a sufficient condition for constructing such matrices is that the SDS $S_1, S_2, S_3, S_4$ with parameters $4 - \{v;k_1, k_2, k_3,k_4;\lambda\}$, satisfy
\begin{displaymath}
    k_1 + k_2 + k_3 + k_4 = v + \lambda
\end{displaymath}
Additionally, if $S_1$ is skew the resulting Hadamard matrix will be skew.

Given these SDS (described in section \ref{sds-computation-section}), we can create $v\times v$ matrices $A_1$, $A_2$, $A_3$, $A_4$, which have entries:
\begin{displaymath}
    (A_n)_{i, j} = 
    \begin{cases}
    -1& \text{if } i-j \in S_n\\
    1              & \text{otherwise}
\end{cases}
\end{displaymath}

These matrices can be plugged into the Goethals-Seidel array to get the desired Hadamard matrix.

\subsubsection*{Skew Hadamard Matrix of Order 292}
I used this construction to obtain a skew Hadamard matrix of order 292. According to multiple papers, this matrix was known since 1978. However, the only paper which we were able to find describing a construction for it was a paper by Djoković \cite{Djokovic2010}, which incorrectly cites the construction for a non-skew Hadamard matrix.

We contacted the author of \cite{Djokovic2010}, who acknowledged that the paper did not contain the correct reference, and helped us by constructing himself the skew SDS that can be used to obtain the skew Hadamard matrix.
\subsection{Construction from T-sequences}
\label{CooperWallisConstruction-section}
Cooper and Wallis \cite{Cooper1972} described a construction of Hadamard matrices that uses ``T-matrices".

\begin{definition}[\cite{Seberry2017}]
    Four circulant $(0,1,-1)$ matrices $X_i$, $i = 1,2,3,4$, of order $n$ which are non-zero for each of the $n^2$ entries for exactly one i, and which satisfy
\begin{displaymath}
    \sum_{i=1}^4X_iX_i^\top = nI
\end{displaymath}
will be called T-matrices of order $n$.
\end{definition}

Given such matrices, and Williamson matrices $A, B, C, D$ of order $w$, Cooper and Wallis showed that it is possible to construct a Hadamard matrix of order $4nw$. Let
\begin{align*}
    e_1 &= GS(X_1, X_2, X_3, X_4) \\
    e_2 &= GS(X_2, -X_1, X_4, -X_3) \\
    e_3 &= GS(X_3, -X_4, -X_1, X_2) \\
    e_4 &= GS(X_4, X_3, -X_2, -X_1)
\end{align*}

Then, the Hadamard matrix is given by ($a \times b$ represent the tensor product between $a$ and $b$):
\begin{displaymath}
    H = e_1 \times A + e_2 \times B + e_3 \times C + e_4 \times D
\end{displaymath}

Although some T-matrices have been computed directly (e.g. in \cite{Sawade1985}), a more convenient way to obtain them is from T-sequences. In fact, if we use T-sequences of length $n$ as the first rows of four circulant matrices $X_1, X_2, X_3, X_4$, we obtain four T-matrices.

Hence, we can look for a Hadamard matrix of order $4n$ by checking if there is a decomposition of $n$ into two factors $w, t$ such that we have Williamson matrices of order $w$ and T-matrices of order $t$ (see Fig. \ref{cooper-wallis-code}).

\begin{lstlisting}[caption={Cooper-Wallis construction}, label=cooper-wallis-code, float]
def hadamard_matrix_cooper_wallis(n):
    for T_seq_len in divisors(n//4):
        will_size = n // (4*T_seq_len)
        if get_T_sequences(T_seq_len, existence=True) and         
          get_williamson_matrices(will_size, existence=True):
          
            x1, x2, x3, x4 = get_T_sequences(T_seq_len)
            a, b, c, d = get_williamson_matrices(will_size)
            
            M = cooper_wallis_construction(x1, x2, x3, x4, 
                                           a, b, c d)
            return M
    return None
\end{lstlisting}

\subsection{Skew Hadamard Matrices from Good Matrices}
A different type of matrices, called Good matrices, can be used to construct skew Hadamard matrices. They are defined in \cite{Koukouvinos2008} as:
\begin{definition}
    Four $(1, -1)$ matrices $A, B, C, D$ of order $n$ (odd) with the properties:
\begin{itemize}
    \item $MN^\top = NM^\top$ for $M, N \in \{A, B, C, D\}$;
    \item $(A-I)^\top = -(A-I)$, $B^\top = B$, $C^\top = C$, $D^\top = D$; 
    \item $AA^\top + BB^\top + CC^\top + DD^\top = 4nI$.
\end{itemize}
will be called good matrices.
\end{definition}

Good matrices of order $n=1, 3, ..., 31$ are listed in \cite{Szekeres1987}. These matrices can be used to obtain a skew Hadamard matrix of order $4n$:

\begin{displaymath}
H = \begin{pmatrix}
A & B & C & D \\
-B & A & D & -C \\
-C & -D & A & B \\
-D & C & -B & A
\end{pmatrix}
\end{displaymath}

\subsection{Miyamoto Construction}
In \cite{Miyamoto1991}, Miyamoto provided a construction for an infinite series of Hadamard matrices:
\begin{theorem}
    Let $q$ be a prime power and $q \equiv 1 \mod 4$. If there is a Hadamard matrix of order $q-1$, then there is a Hadamard matrix of order $4q$. \qed
\end{theorem}

This construction uses conference matrices, which are matrices where the diagonal contains all zeros, all other entries are $(-1, 1)$ and the matrix satisfies $CC^\top = (n-1)I$.

Paley's second construction of Hadamard matrices provides a method to obtain conference matrices of order $q+1$ ($q$ an odd prime power) in the form:
\begin{displaymath}
C = \begin{pmatrix}
0 & e_{q} \\
e_{q}^\top & D
\end{pmatrix}
\end{displaymath}
where $e_k$ is the row vector containing all ones. 

Then, by permuting rows and columns of the matrix, we can rewrite it as:
\begin{displaymath}
    C = \begin{pmatrix}
    0 &  1 & e_{m} & e_m \\
    1 &  0 & e_{m} & -e_m \\
    e_m^\top &  e_m^\top & -C_1 & -C_2 \\
    e_m^\top &  -e_m^\top & C_2^\top & C_4
    \end{pmatrix}
\end{displaymath}

Now, the Hadamard matrix $K$ of order $q-1$ is split into four sub matrices $K_1$, $K_2$, $K_3$, $K_4$:
\begin{displaymath}
K = \begin{pmatrix}
K_1 & K_2 \\
-K_3 & K_4
\end{pmatrix}
\end{displaymath}

To simplify notation, define the following:
\begin{align*}
U_{11} &= U_{33} = C_1 \\
U_{12} &= U_{34} = C_2 \\
U_{21} &= U_{43} = -C_2^\top \\
U_{22} &= U_{44} = C_4 \\
U_{13} &= U_{14} = U_{23} = U_{24} = U_{31} = U_{32} = U_{41} = U_{42} = 0 \\
V_{13} &= -V_{31}^\top = K_1 \\
V_{14} &= -V_{41}^\top = K_2 \\
V_{23} &= -V_{32}^\top = K_3 \\
V_{24} &= V_{32}^\top = K_4 \\
V_{11} &= V_{22} = V_{33} = V_{44} = I \\
V_{12} &= V_{21} = V_{34} = V_{43} = 0 \\
\end{align*}

and then, for $1 \le i,j \le 4$ let:

\begin{displaymath}
T_{ij} = \begin{pmatrix} 
U_{ij} + V_{ij}& U_{ij} - V_{ij} \\
U_{ij} - V_{ij} & U_{ij} + V_{ij}
\end{pmatrix}
\end{displaymath}

Finally, the Hadamard matrix of order $4q$ is given by:
\begin{displaymath}
H = \begin{pmatrix}
1 & -e & 1 & e & 1 & e & 1 & e \\
-e^\top & T_{11} & e^\top & T_{12} & e^\top & T_{13} & e^\top & T_{14} \\
-1 & -e & 1 & -e & 1 & e & -1 & -e \\
-e^\top & -T_{21} & -e^\top & T_{22} & e^\top & T_{23} & -e^\top & -T_{24} \\
-1 & -e & -1 & -e & 1 & -e & 1 & e \\
-e^\top & -T_{31} & -e^\top & -T_{32} & -e^\top & T_{33} & e^\top & T_{34} \\
-1 & -e & 1 & e & -1 & -e & 1 & -e \\
-e^\top & -T_{41} & e^\top & T_{42} & -e^\top & -T_{43} & -e^\top & T_{44}
\end{pmatrix}
\end{displaymath}
\subsection{Spence Construction from Supplementary Difference Sets}
Spence \cite{Spence1975} proposed a construction that uses SDS with parameters $4-\{2v; v, v, v, v+1; 2v\}$:
\begin{theorem}
    If there exist $4-\{2v; v, v, v, v+1; 2v\}$ supplementary difference sets in the cyclic group of order $2v$ then there exists a Hadamard matrix of order $4(2v+1)$. \qed
\end{theorem}

The SDS $S_1, S_2, S_3, S_4$ are used to create four matrices $A_1, A_2, A_3, A_4$:
\begin{displaymath}
    (A_l)_{i,j} = \begin{cases}
    +1& \text{if } i-j \in S_l\\
    -1&  \text{if } i-j \not\in S_l
\end{cases}
\end{displaymath}

Furthermore, let $P$ be the permutation matrix with ones on the anti-diagonal. The Hadamard matrix is given by:
\begin{displaymath}
H = \begin{pmatrix}
+1 & -1 & +1 & +1 &  e &  e &  e &  e \\
+1 & +1 & -1 & +1 & -e &  e & -e &  e \\
-1 & +1 & +1 & +1 & -e &  e &  e & -e \\
-1 & -1 & -1 & +1 & -e & -e &  e &  e \\
-e^\top &  e^\top &  e^\top & -e^\top &  A_1  &  A_2P & A_3P & A_4P \\
-e^\top & -e^\top &  e^\top &  e^\top & -A_2P &  A_1  & -A_4^\top P & A_3^\top P \\
-e^\top & -e^\top & -e^\top & -e^\top & -A_3P &  A_4^\top P & A_1 & -A_2^\top P \\
 e^\top & -e^\top &  e^\top & -e^\top & -A_4P & -A_3^\top P & A_2^\top P & A_1 
\end{pmatrix}
\end{displaymath}
\subsection{Construction from Complementary Difference Sets}
In \cite{Blatt1969}, Blatt and Szekeres showed that given complementary difference sets of size $m$, a skew Hadamard matrix of order $4(m+1)$ exists. 

Suppose that $A, B$ are two complementary difference sets over a group $G$. Then, the Hadamard matrix will be in the form $H = I+S$, where $S$ is defined as follows. Let  $\gamma_1, \gamma_2, ..., \gamma_{2m+1}$ be the elements of $G$. We have, for $1 \le i,j \le 2m+1$:
\begin{align*}
    -S_{2m+1+i, 2m+1+j} = S_{ij}  =  \begin{cases}
                            +1& \text{if } \gamma_j-\gamma_i \in A\\
                            -1 & \text{otherwise}
                        \end{cases} \\
    -S_{2m+1+j, i} = S_{i, 2m+1+j}  =  \begin{cases}
                            +1& \text{if } \gamma_j-\gamma_i \in B\\
                            -1              & \text{otherwise}
                        \end{cases} \\                   
\end{align*}

Furthermore, 

\begin{align*}
    -S_{4m+3, i} &= S_{i,4m+3} =  \begin{cases}
                            +1& \text{for }1 \le i\le 2m+1\\
                            -1& \text{for }2m+2 \le i\le 4m+2
                        \end{cases} \\
    S_{4m+4, i} &= -S_{i, 4m+4} = +1 \text { for } 1\le i \le 4m+3 \\
    S_{ii} &= 0 \text{ for } 1 \le i \le 4m+4
\end{align*}

It is easy to see from the definition of $S$ that it is skew-symmetric. Hence, the resulting Hadamard matrix will be skew.
\subsection{Spence Construction of Skew Hadamard Matrices}
Another construction of skew Hadamard matrices using complementary difference sets was given by Spence \cite{Spence1975b}.

\begin{theorem}
    If there exists a cyclic projective plane of order $q$ and two complementary difference sets in a cyclic group of order $1 + q + q^2$, then there exists a skew-Hadamard matrix of the Goethals-Seidel type of order $4(1 + q + q^2)$. \qed
\end{theorem}

He noted that cyclic projective planes of order $q$ always exist if $q$ is a prime power. Therefore, using the known constructions of complementary difference sets it is easy to see that a skew Hadamard matrix of order $4(1+q+q^2)$ can be constructed whenever $q$ is a prime power such that either $1+q+q^2$ is a prime congruent to $3, 5, 7 \mod 8$ or $2q^2+2q+3$ is a prime power.

The cyclic projective plane is used by Spence to create a $(1 + q^2 + q^4, 1 + q^2, 1)$ difference set. Since a construction for such sets Was already present in SageMath, I used it directly. Then, let $D$ be a translate of this set fixed by $q$. We define a subset $D_1$:
\begin{displaymath}
    D_1 = \{d \in D \text{ } | \text{ } d \equiv 0 \mod 1 - q + q^2 \}
\end{displaymath}

The elements in $D \setminus D_1$ can be partitioned into pairs $(d_i, d'_i)$ such that $d_i \equiv d'_i \mod 1+q+q^2$ and  $d_i \not \equiv d_j \mod 1+q+q^2$  whenever $i \not = j$. Given these pairs, we define a new set:
\begin{displaymath}
    D_2 = \{ d_i \mod 1+q+q^2 \text{ } | 1 \le i \le \frac{1}{2}(1+q^2)\}
\end{displaymath}

Then, we obtain the matrices $R$ and $S$ of order $1+q+q^2$:
\begin{align*}
    R_{ij} &= \begin{cases}
                +1& \text{if } j-i \equiv d \mod 1+q+q^2 \text{ for some } d\in D_1 \cup D_2\\
                -1 & \text{otherwise}
            \end{cases} \\
    S_{ij} &= \begin{cases}
                +1& \text{if } j-i \equiv d \mod 1+q+q^2 \text{ for some } d\in D_2\\
                -1 & \text{otherwise}
            \end{cases} \\
\end{align*}

Now, let $A, B$ be the complementary difference sets of order $1+q+q^2$. We define two additional matrices of order $1+q+q^2$:

\begin{align*}
    P_{ij} &= \begin{cases}
                +1& \text{if } j-i \in A\\
                -1 & \text{otherwise}
            \end{cases} \\
    Q_{ij} &= \begin{cases}
                +1& \text{if } j-i \in B\\
                -1 & \text{otherwise}
            \end{cases} \\
\end{align*}

Finally, plugging the matrices $P, Q, R, S$ into the Goethals-Seidel array will give a skew Hadamard matrix of order $4(1+q+q^2)$.
\subsection{Skew Hadamard Matrices from Amicable Orthogonal Designs}
The last construction of skew Hadamard matrices that I implemented uses the following theorem \cite{Seberry1978}:

\begin{theorem}
    Suppose there is an orthogonal design of type $(1,m,mn-m-1)$ in order $mn$. Suppose $n$ is the order of amicable designs of types $((1, n-1); (n))$. Then there is a skew Hadamard matrix of order $mn(n-1)$. \qed
\end{theorem}

Orthogonal designs are defined as follows:
\begin{definition}[Seberry \cite{Seberry1978}]
    An orthogonal design of order $n$ and type $(u_1, u_2, ..., u_k)\in \mathbb{Z}^k_{>0}$  on the commuting variables $x_1, x_2, ..., x_k$ is an $n \times n$ matrix A with entries from $\{0, x_1, x_2, ...,x_k\}$ such that
    \begin{displaymath}
        AA^\top = \sum_{i=1}^k(u_ix_i^2)I_n
    \end{displaymath}\qed
\end{definition}

Furthermore, two orthogonal designs $M$ of type $(m_1, ..., m_p)$ and $N$ of type $(n_1, .., n_q)$, both  of order $n$, are called amicable of type $((m_1, ...,m_p); (n_1, ..., n_q))$ if $MN^\top = NM^\top$.

Now, suppose $M$, $N$ are amicable orthogonal designs of type $((1, n-1); (n))$. They can be rewritten in the form:
\begin{equation}
\label{amicable-orth-design-normal-form}
    M = xI + y\begin{pmatrix}
        0 & e  \\
        -e^\top & P
    \end{pmatrix}
    \quad \quad
    N = z\begin{pmatrix}
        1 & e  \\
        -e^\top & D
    \end{pmatrix}
\end{equation}
where $e$ is the $1\times n-1$ vector of all ones. 

Then, let $J$ be the $(n-1) \times (n-1)$ matrix of all ones, and $K$ be the orthogonal design of order $mn$ and type $(1,m,mn-m-1)$. Replacing the variables of $K$ with $P$, $J$, $D$ respectively we get a skew Hadamard matrix of order $mn(n-1)$.

I used this construction to obtain a skew Hadamard matrix of order $756$, by setting $m=1$ and $n=28$. However, to do so it has been necessary to find the orthogonal designs $M$, $N$, $K$ used in the construction. 

 An orthogonal design of type $(1, 1, 26)$ ($K$) is listed in Appendix G of \cite{Seberry2017}. However, for the amicable orthogonal designs $M$, $N$ of type $((1, 27);(28))$ I used a more general approach, which creates orthogonal designs of type $((1, n-1); (n))$ from amicable Hadamard matrices of order $n$.

\subsubsection{Construction of Amicable Hadamard Matrices}
Amicable Hadamard matrices are defined as:
\begin{definition}[Seberry\cite{Seberry2017}, Chapter 5.10]
    Two $n\times n$ Hadamard matrices $W$ and $M$ are called amicable Hadamard matrices of order $n$ if:
    \begin{itemize}
        \item $W$ is a skew Hadamard matrix;
        \item $M$ is symmetric;
        \item $WM^\top = MW^\top$.
    \end{itemize}
\end{definition}

Seberry \cite{Seberry2017} proved that two such matrices $M, W$ of order $n$ are equivalent to two amicable orthogonal designs of type $((1, n-1);(n))$. If we write $W = I + S$, the two orthogonal designs will be:
\begin{displaymath}
    A = xI+yS \quad \quad B = zM
\end{displaymath}

Furthermore, if $M$ and $W$ are normalised, the resulting orthogonal designs will be in the form described in Eq. \ref{amicable-orth-design-normal-form}.

Hence, we have now transformed the problem into finding two amicable Hadamard matrices of order 28. A construction of such matrices is described by Seberry in \cite{Wallis1970}: this works for any order $n=q+1$, where $q$ is a prime power.

Let $G$ be a Galois field of order $q$, and order its elements such that $a_0=0$ and $a_{q-i} = -a_i$ for $1 \le i \le q-1$. Then, define $\chi(x)$ as follows:
\begin{displaymath}
    \chi(x) = \begin{cases}
                0& \text{if } x=0 \\
                1& \text{if } x \text{ is a square} \\
                -1& \text{otherwise}
            \end{cases}
\end{displaymath}

We can create two $q \times q$ matrices $S, R$, where $S$ is defined by $S_{ij} = \chi(a_j-a_i)$, and the entries of $R$ are:
\begin{displaymath}
        R_{ij} = \begin{cases}
        1& \text{if } i=j=0 \lor (j=q-i \land 1 \le i \le q-1) \\
        0& \text {otherwise}
    \end{cases}
\end{displaymath}

Now, let $P = S+I$ and $D = R + RS$. The two amicable Hadamard matrices are given by:
\begin{displaymath}
    W = \begin{pmatrix}
        1 & e \\
        -e^\top & P
    \end{pmatrix}
    \quad \quad
    M = \begin{pmatrix}
        1 & e \\
        e^\top & D
    \end{pmatrix}
\end{displaymath}

\clearpage
\section{Conclusion}
\label{conclusion-chapter}
The aim of this project was to implement Hadamard and skew Hadamard matrices of order up to 1000 in SageMath. I achieved this by implementing nine new methods for obtaining these matrices (as well as some additional constructions, as seen in chapters \ref{NA-chapter} and \ref{diff-sets-chapter}). Therefore,  from version 10.0 of SageMath it will be possible to construct skew and non-skew Hadamard matrices of every order less than or equal to 1000, for which a construction is known.

In total, I have added more than 5000 lines of code to the SageMath repository (of which around 1000 contain the documentation). This code is easily extendable: as  explained in section \ref{utility-funcs-section}, a generic  function {\tt hadamard\_matrix} which accesses all implemented constructions is available. Therefore, if a new construction is added to SageMath,  it is enough to add it also to this function, and then all users that were using Hadamard matrices will be able to use the new orders without having to modify their code.

Furthermore, extensive tests have been written for all the methods added to SageMath, and users can use the {\tt check} parameter if they want to be sure that the function is working correctly. These design decisions make sure that the requirements detailed in Chapter \ref{intro-chapter} have been met.

\subsection{Challenges}
While working on this project, one of the main difficulties has been that many of the papers that I read use different notations, often conflicting with each other. For example, the name ``complementary difference sets" has sometimes been used to refer to supplementary difference sets. 

Furthermore, it has been hard to determine how to create some Hadamard matrices, because  many papers listed the known orders without citing the relevant constructions. In particular, this happened with the skew Hadamard matrix of order 292, for which we could not find any published construction.

Lastly, while implementing some functions I  realised that they were generalisations of constructions that I had already added to SageMath. Therefore, if I had to work on this project again I would try to choose the constructions more carefully, to avoid redundancy.

\subsection{Future Work}
Given the nature of the project,  it would be easy to extend it by adding more constructions for Hadamard matrices. On the other hand, it may also be useful to search for constructions that cover new orders of Hadamard matrices.

Lastly, since the Hadamard conjecture is still open, it would be interesting to try to prove (or disprove) it.

\appendix
\section{Table Of Constructions}
Table \ref{hadamard-constructions-list-table} details for every odd value of $n < 250$ which method can be used to construct the Hadamard matrix of order $4n$ (whenever $n$ is even, the Hadamard matrix can be constructed using the doubling construction). Similarly, Table \ref{skew-hadamard-constructions-list-table} lists the constructions of skew Hadamard matrices. Note that some entries of the two tables are empty: this indicates that no construction for the corresponding order is known.
Table \ref{construction-list-table} contains an explanation of the abbreviations used in the two lists. 

\renewcommand{\arraystretch}{1.5}
\begin{table}[hbt]
\caption{Algorithms used in tables 2-3}
\label{construction-list-table}
\begin{centering}
\begin{tabular}{>{\raggedleft\arraybackslash}p{0.20\linewidth} | >{\raggedright\arraybackslash}p{0.75\linewidth}}
 PaleyI &  Paley's first construction \\
 PaleyII & Paley's second construction\\  
 Will &  Williamson construction \\
 GS & Goethals-Seidel array \\
 SDS & Construction from (possibly skew) SDS \\
 CW(t) & Construction of Hadamard matrix of order $4n$ from T-sequences of length $t$ (and Williamson matrices of order $n/t$) \\
 Good & construction from good matrices \\
 Miy & Miyamoto Construction \\
 CDS & Construction from complementary difference sets \\
 Spence(q) & Spence construction of skew Hadamard matrix  of order $4(1+q+q^2)$ \\
 AOD(m, n) & construction of a skew Hadamard matrix of order $mn(n-1)$ from amicable orthogonal designs \\
 \end{tabular}
 \end{centering}    
\end{table}

\begin{table}[hbt]
    \caption{Hadamard matrices of order $4n$, up to 1000}
    \label{hadamard-constructions-list-table}
    \centering
    \begin{tabular}{|c c || c c || c c || c c || c c|}
        1 & PaleyI &
        3 & PaleyI &
        5 & PaleyI &
        7 & PaleyI &
        9 & PaleyII \\
        11 & PaleyI &
        13 & PaleyII &
        15 & PaleyI &
        17 & PaleyI &
        19 & PaleyII \\
        21 & PaleyI &
        23 & Will &
        25 & PaleyII &
        27 & PaleyI &
        29 & Will \\
        31 & PaleyII &
        33 & PaleyI &
        35 & PaleyI &
        37 & PaleyII &
        39 & Will \\
        41 & PaleyI &
        43 & Will &
        45 & PaleyI &
        47 & CW(47) &
        49 & PaleyII \\
        51 & PaleyII &
        53 & PaleyI &
        55 & PaleyII &
        57 & PaleyI &
        59 & CW(59) \\
        61 & PaleyI &
        63 & PaleyI &
        65 & CW(5) &
        67 & CW(67) &
        69 & PaleyII \\
        71 & PaleyI &
        73 & Miy &
        75 & PaleyII &
        77 & PaleyI &
        79 & PaleyII \\
        81 & CW(3) &
        83 & PaleyI &
        85 & PaleyII &
        87 & PaleyI &
        89 & CW(89) \\
        91 & PaleyII &
        93 & CW(3) &
        95 & PaleyI &
        97 & PaleyII &
        99 & PaleyII \\
        101 & Miy &
        103 & SDS &
        105 & PaleyI &
        107 & CW(107) &
        109 & Miy \\
        111 & PaleyI &
        113 & Miy &
        115 & PaleyII &
        117 & PaleyI &
        119 & CW(7) \\
        121 & PaleyII &
        123 & PaleyI &
        125 & PaleyI &
        127 & SDS &
        129 & PaleyII \\
        131 & PaleyI &
        133 & CW(7) &
        135 & PaleyII &
        137 & PaleyI &
        139 & PaleyII \\
        141 & PaleyI &
        143 & PaleyI &
        145 & PaleyII &
        147 & PaleyI &
        149 & Miy \\
        151 & SDS &
        153 & CW(3) &
        155 & PaleyI &
        157 & PaleyII &
        159 & PaleyII \\
        161 & PaleyI &
        163 & SDS &
        165 & PaleyI &
        167 &  &
        169 & PaleyII \\
        171 & PaleyI &
        173 & PaleyI &
        175 & PaleyII &
        177 & PaleyII &
        179 &  \\
        181 & PaleyII &
        183 & CW(3) &
        185 & PaleyI &
        187 & PaleyII &
        189 & CW(3) \\
        191 & SDS &
        193 & Miy &
        195 & PaleyII &
        197 & PaleyI &
        199 & PaleyII \\
        201 & PaleyII &
        203 & PaleyI &
        205 & PaleyII &
        207 & PaleyI &
        209 & CW(11) \\
        211 & PaleyII &
        213 & CW(71) &
        215 & PaleyI &
        217 & PaleyII &
        219 & SDS \\
        221 & PaleyI &
        223 &  &
        225 & PaleyII &
        227 & PaleyI &
        229 & PaleyII \\
        231 & PaleyII &
        233 & Miy &
        235 & CW(47) &
        237 & PaleyI &
        239 & SDS \\
        241 & Miy &
        243 & PaleyI &
        245 & CW(5) &
        247 & CW(13) &
        249 & CW(83) 
    \end{tabular}
\end{table}

\begin{table}[hbt]
    \caption{Skew Hadamard matrices of order $4n$, up to 1000}
    \label{skew-hadamard-constructions-list-table}
    \centering
    \begin{tabular}{|c c || c c || c c || c c || c c|}
        1 & Good &
        3 & Good &
        5 & Good &
        7 & Good &
        9 & Good \\
        11 & Good &
        13 & Good &
        15 & Good &
        17 & Good &
        19 & Good \\
        21 & Good &
        23 & Good &
        25 & Good &
        27 & Good &
        29 & Good \\
        31 & Good &
        33 & PaleyI &
        35 & PaleyI &
        37 & SDS &
        39 & SDS \\
        41 & PaleyI &
        43 & SDS &
        45 & PaleyI &
        47 & GS &
        49 & SDS \\
        51 & CDS &
        53 & PaleyI &
        55 & CDS &
        57 & PaleyI &
        59 & GS \\
        61 & PaleyI &
        63 & PaleyI &
        65 & SDS &
        67 & SDS &
        69 & GS \\
        71 & PaleyI &
        73 & SDS &
        75 & CDS &
        77 & PaleyI &
        79 & CDS \\
        81 & SDS &
        83 & PaleyI &
        85 & CDS &
        87 & PaleyI &
        89 &   \\
        91 & CDS &
        93 & SDS &
        95 & PaleyI &
        97 & SDS &
        99 & CDS \\
        101 &   &
        103 & SDS &
        105 & PaleyI &
        107 &   &
        109 & SDS \\
        111 & PaleyI &
        113 & SDS &
        115 & CDS &
        117 & PaleyI &
        119 &   \\
        121 & SDS &
        123 & PaleyI &
        125 & PaleyI &
        127 & SDS &
        129 & SDS \\
        131 & PaleyI &
        133 & SDS &
        135 & CDS &
        137 & PaleyI &
        139 & CDS \\
        141 & PaleyI &
        143 & PaleyI &
        145 & SDS &
        147 & PaleyI &
        149 &   \\
        151 & SDS &
        153 &   &
        155 & PaleyI &
        157 & SDS &
        159 & CDS \\
        161 & PaleyI &
        163 & SDS &
        165 & PaleyI &
        167 &   &
        169 & SDS \\
        171 & PaleyI &
        173 & PaleyI &
        175 & CDS &
        177 &   &
        179 &   \\
        181 & SDS &
        183 & Spence(13) &
        185 & PaleyI &
        187 & CDS &
        189 & AOD(1, 28) \\
        191 &   &
        193 &   &
        195 & CDS &
        197 & PaleyI &
        199 & CDS \\
        201 &   &
        203 & PaleyI &
        205 &   &
        207 & PaleyI &
        209 &   \\
        211 & CDS &
        213 & SDS &
        215 & PaleyI &
        217 & SDS &
        219 & SDS \\
        221 & PaleyI &
        223 &   &
        225 &   &
        227 & PaleyI &
        229 &   \\
        231 & CDS &
        233 &   &
        235 &   &
        237 & PaleyI &
        239 & SDS \\
        241 & SDS &
        243 & PaleyI &
        245 &   &
        247 & SDS &
        249 &   
    \end{tabular}
\end{table}

\clearpage

\printbibliography[heading=bibintoc]
\end{document}